\documentclass{amsart}
  \usepackage{amsmath,amssymb,amsthm,amsfonts,graphicx,color}
  \definecolor{limegreen}{rgb}{0.196,0.804,0.196}
  \definecolor{darkgreen}{rgb}{0.0,0.5,0.0}
  \definecolor{darkbluegreen}{rgb}{0,0.3,0.6}
  \definecolor{badgerred}{rgb}{0.715,0.004,0.004}
  \usepackage[font=small,labelfont={bf,sf}, textfont={sf},margin=0em]{caption}
  \usepackage{url,hyperref}
  \hypersetup{colorlinks, 
              citecolor=badgerred,
              filecolor=black,
              linkcolor=blue,
              urlcolor=darkgreen,
              bookmarksopen=false}

  \newcommand{\vN}{{\boldsymbol N}}
  
  \newcommand{\vT}{{\boldsymbol T}}

  \newcommand{\cA}{\mathcal{A}}
  \newcommand{\cB}{\mathcal{B}}
  \newcommand{\cC}{\mathcal{C}}
  \newcommand{\Ca}{\Cb}
  \newcommand{\Cb}{\bar{\mathcal{C}}}

  \newcommand{\cG}{\mathcal{G}}

  \newcommand{\cO}{\mathcal{O}}
  
  \newcommand{\dd}{\mathrm{d}}
  \newcommand{\pd}{\partial}
  \newcommand{\R}{{\mathbb R}}

  \title[Ancient Solutions to Curve Shortening]{Ancient Solutions to Curve Shortening\\with Finite Total Curvature}
  \author{Sigurd Angenent}
  \address{Department of Mathematics, UW Madison}
  \author{Qian You}
  \address{Jersey City, NJ07310}
  \begin{document}

\begin{abstract}
  We construct ancient solutions to Curve Shortening in the plane whose total curvature is uniformly bounded by gluing together an arbitrary chain of given Grim Reapers along their common asymptotes.
\end{abstract}

\maketitle
\section{Introduction}

\subsection{Ancient solutions}
A family of plane curves $\{\cC_t : t_0<t<t_1\}$ moves by Curve Shortening if for some parametrization $\cC:\R\times(t_0, t_1) \to \R^2$ one has
\begin{equation} \label{eq-CS} \left(\frac{\pd\cC}{\pd t}\right)^\perp = \frac{\pd^2 \cC}{\pd s^2}.
\end{equation}
Here $(\pd_t\cC)^\perp$ is the component of the velocity vector $\cC_t$ which is perpendicular to the curve, and $\pd/\pd s$ stands for the arc length derivative along the curve.

An ancient solution of Curve Shortening is a solution $\cC_t$ that is defined for all $t\in(-\infty, t_1)$, for some $t_1\in\R$. While it can be shown that for any reasonably smooth initial plane curve $\cC_0$ a unique solution $\{\cC_t : 0<t<t_1\}$ to Curve Shortening exists with $\cC_0$ as initial curve, the requirement that a solution be defined for all times $t\in(-\infty, t_1)$ for some given $t_1$ is much more restrictive. Two ancient solutions that have been known for a long time are the shrinking circle (with radius $r(t) = \sqrt{2(T-t)}$, and some arbitrary but fixed center), and the translating soliton known as the \textit{Grim Reaper}. This last curve is up to rotation and translation defined by the equation
\begin{equation} \label{eq-grim-reaper}
  x-C = \frac{1}{v}G\bigl(v(y-a)\bigr) + vt, \text{ for } 0<y<\frac\pi v, t\in\R,
\end{equation}
where $v>0$ is the velocity of the soliton, and where, by definition,
\[
  G(y) = -\ln\sin y.
\]

\begin{figure}
  \centering
  \includegraphics[width=0.6\textwidth]{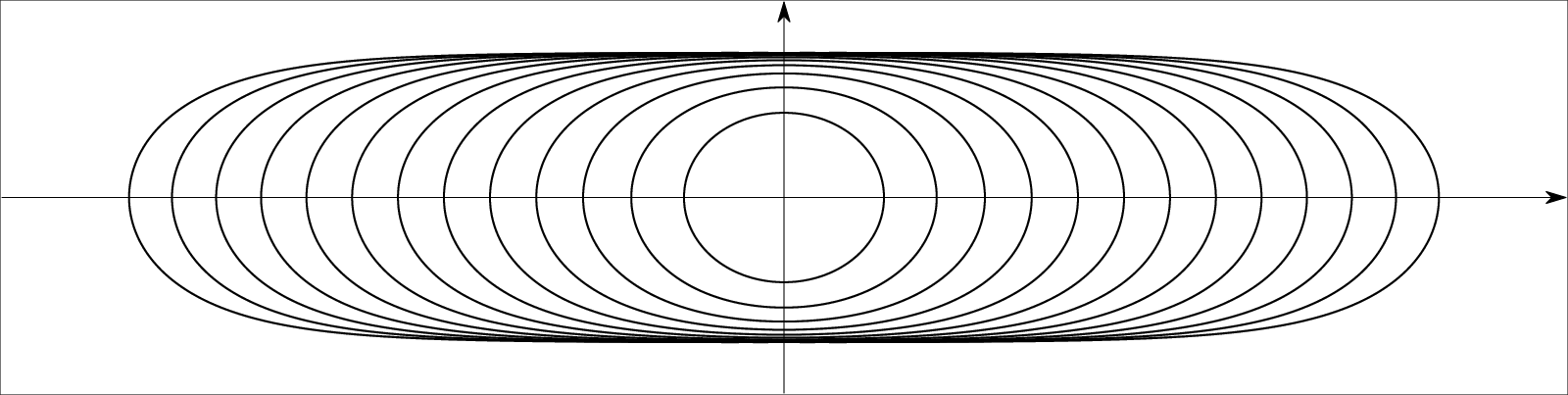}\\[2ex]
  \includegraphics[width=0.6\textwidth]{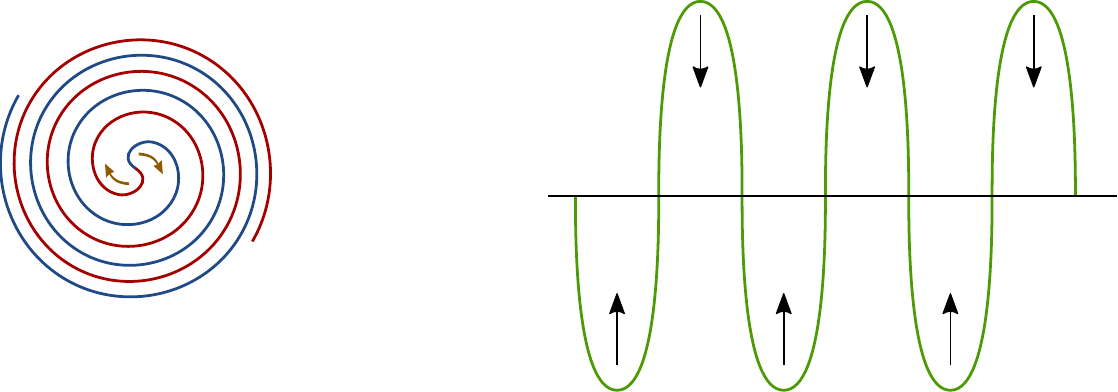}
  \caption{Ancient convex solutions to Curve Shortening: the paperclip solution (top), the Yin-Yang rotating soliton \cite{Altschuler93} (bottom left), and the Ancient Sine Curve \cite{nakayamalizukawadati:cle94} (bottom right).}
  \label{fig-ancient-convex}
\end{figure}

There is a short list of other known examples of ancient solutions to Curve Shortening. The \emph{Abresch-Langer curves} are compact, immersed curves that shrink self similarly to their center of mass. Except for the circle, none of these curves are embedded. The \emph{paper-clip solution,} found by Angenent \cite{angenent:shrinkingdoughnuts} and also by Nakayama, Iizuka, and Wadati~\cite{nakayamalizukawadati:cle94} (Figure~\ref{fig-ancient-convex}, top) is a convex embedded curve which for $t\to-\infty$ is asymptotically described by two Grim Reaper solutions with the same asymptotes, moving toward each other. Another ancient solution is the \emph{Yin-Yang curve,} which is a spiral shaped so that it evolves simply by rotating at a steady rate. (The curve was introduced by Altschuler in \cite{Altschuler93}.  See also \cite{2000HungerBuehlerSmoczyk,MR3043378}.) There is also the \emph{Ancient Sine Wave} discovered by Nakayama, Iizuka, and Wadati (see \cite{nakayamalizukawadati:cle94} and also Figure~\ref{fig-ancient-convex}.)  Both the paper-clip and the ancient sine wave have explicit parametrizations.

One can try to order ancient solutions by their total curvature.  The total curvature of a plane curve $\cC$ (compact or not) is defined to be
\[
  K(\cC) = \int_{\cC} |\kappa|\, \dd s
\]
It is well known that the total curvature is a monotone quantity along a family of curves $\cC_t$ that evolve by Curve Shortening \cite{grayson:hes86,angenent:pec90}.  One has
\[
  \frac{\dd}{\dd t} K(\cC_t) = - 2 \sum_{P_i(t)} \left|\frac{\pd\kappa}{\pd s}\Bigr|_{P_i}\right|
\]
where the sum is taken over all inflection points $P_i(t)$ of the curve $\cC_t$.

The total curvature of a closed curve is always at least $2\pi$, and a closed curve $\cC$ is convex if and only if its total curvature is exactly $2\pi$.  The Daskalopoulos-Hamilton-Sesum theorem \cite{DaskalopoulosHamiltonSesum:CompactAncientCS} is a classification of all ancient solutions of this type.  They showed that any ancient embedded, convex, and compact solution must be either a shrinking circle, or else the ancient paper-clip solution.

At the other extreme there are ancient solutions of infinite total curvature, e.g.~the Yin-Yang curve, and the Ancient Sine Wave.

\subsection{Main result}
\label{sec-main-result}
We construct a large number of ancient solutions with uniformly bounded total curvature.  To be specific, let two sets of numbers $\{a_0, \dots, a_n\}$ and $\{C_1, \dots, C_{n}\}$ be given.  Then we will construct an ancient solution $\{\cC_t \mid t<0\}$ of Curve Shortening that crosses the $y$-axis at $n+1$ different points $\{P_0(t), \dots, P_n(t)\}$ and for which for each $k=1, 2, \dots, n$, the arc $P_{k-1}(t)P_k(t)$ converges to the translating soliton $\cG_k(t)$ defined by
\begin{equation} \label{eq-asymptotic-to-GR} \cG_k(t) :\qquad x - C_k = G_k(y, t) = (-1)^k \left\{ \frac{1}{v_k} G\bigl(v_k |y-a_{k-1}|\bigr) + v_kt \right\}
\end{equation}
as $t\to-\infty$.  Here
\[
  v_k = \frac{\pi}{|a_k - a_{k-1}|}
\]
is the asymptotic velocity of the $k^{\rm th}$ arc in the solution $\cC_t$ that we construct.

\begin{figure}
  \includegraphics{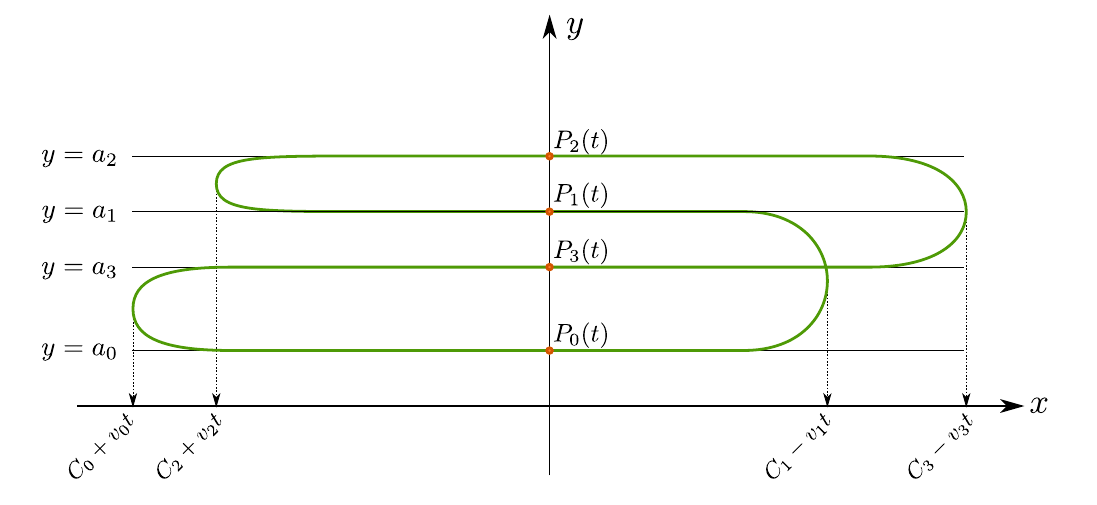}
  \includegraphics{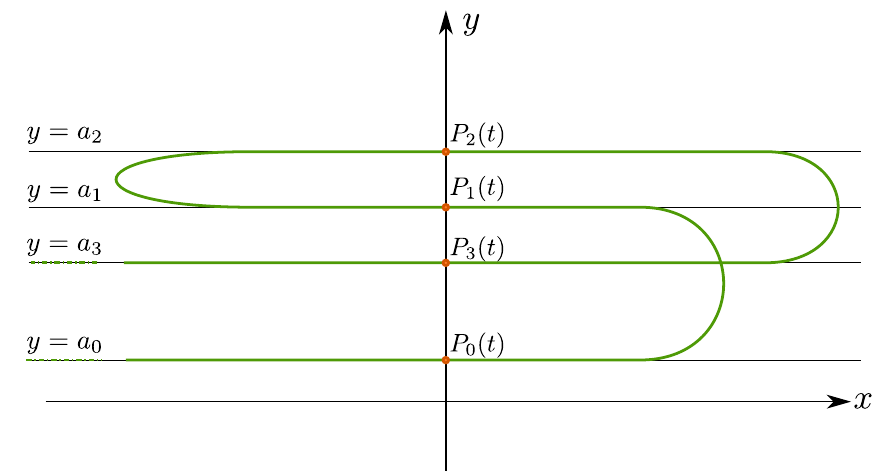}
  \caption{Two ancient solutions constructed in this paper from the same data $\{a_0, a_1, a_2, a_3\}$ and $\{C_0, C_1, C_2, C_3\}$. The solution on top is compact, the bottom solution is asymptotic to the two horizontal lines $y=a_0$, $y=a_n$.}
  \label{fig-TwoDoodles}
\end{figure}

The solutions we construct come in two varieties: compact, and non compact.  To obtain compact solutions we assume that there is a final arc $P_n(t)P_0(t)$ connecting the last and first intersection points.  As $t\to-\infty$ this arc is asymptotic to the translating soliton $x=G_{n+1}(y, t)$, given by \eqref{eq-asymptotic-to-GR}, with velocity $v_0 = v_{n+1} = \pi/|a_0-a_{n}|$.

To obtain non-compact solutions we let the initial arc ending at $P_0(t)$ be asymptotic to the straight line $y=a_0$, and similarly, we let the final arc be asymptotic to the line $y=a_n$.  See Figure~\ref{fig-TwoDoodles}.

For any compact immersed curve the number of intersections with a straight line, such as the $y$-axis, must be even.  Therefore one can only construct non-compact solutions when $n$ is odd.

We will construct these solutions in three stages.  It turns out that convex ancient solutions are easier to construct.  An immersed curve in the plane is convex if it has no inflection points, i.e.~if the curvature never vanishes on the curve.  Such curves need not be embedded (e.g.~consider a cardioid).  Thus we first construct ancient convex immersed solutions in section~\ref{sec-ancient-convex}.  Then, in section~\ref{sec-embedded-ancient-non-convex}, we construct embedded ancient solutions.  Finally, in section~\ref{sec-general-ancient-non-convex} we combine the results and arguments from sections~\ref{sec-ancient-convex} and~\ref{sec-embedded-ancient-non-convex} to deal with the most general case of ancient solutions that are neither embedded nor convex.

\section{Construction of ancient convex solutions}
\label{sec-ancient-convex}

\subsection{The convexity assumption}
In this section we assume that the asymptotes appear in alternating order, i.e.
\begin{equation}
  \label{eq-convexity-hypothesis}
  \textit{$\forall k\in\{1, \dots, n\}$ either }
  a_{k} < \min \left\{ a_{k-1}, a_{k+1}\right\}
  \textit{ or }
  a_{k} > \max \left\{a_{k-1}, a_{k+1}\right\}
\end{equation}
in which we agree to define $a_{n+1}=a_0$.  If we define Grim Reapers as in \eqref{eq-asymptotic-to-GR} then any two consecutive Grim Reapers will intersect:
\subsection{Lemma}\itshape
\label{lem-intersection-location}
Under the convexity condition \eqref{eq-convexity-hypothesis} there exists a $t_0\in\R$ such that for all $t<t_0$ and for every $k$ there is a unique intersection point $A_k(t)$ whose coordinates $(x,y)$ satisfy
\[
  x = G_k(y, t) = G_{k+1}(y, t)
\]
with $|y-a_k| < \frac12 \min \bigl\{|a_{k-1}-a_k|, |a_{k+1}-a_k|\bigr\}$.

\begin{figure}[h]
  \includegraphics[width=0.6\textwidth]{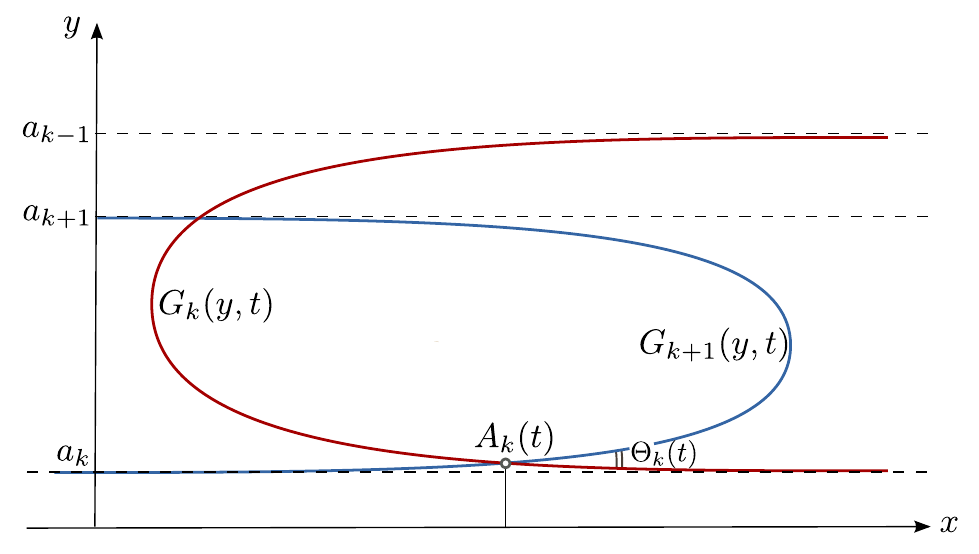}
  \caption{The intersection point $A_k(t)$ and angle $\Theta_k(t)$ of two consecutive Grim Reapers under the convexity hypothesis \eqref{eq-convexity-hypothesis}.}
  \label{fig-intersectionpoint-convex}
\end{figure}
For $t\to-\infty$ the $x$-coordinate of the intersection point $A_k(t)$ satisfies
\begin{equation}
  \label{eq-intersection-convex-x-coord}
  x = \bigl(v_k-v_{k+1}\bigr) t + c_k + \cO\bigl(e^{2v_kv_{k+1}t}\bigr)
\end{equation}
for a suitable constant $c_k$.  \upshape\medskip

\subsection{Lemma}\itshape
\label{lem-intersection-angle}
The angle between the tangents to $\cG_k(t)$ and $\cG_{k+1}(t)$ at the intersection $A_k(t)$ satisfies
\begin{equation}
  \label{eq-intersection-angle-asymptotics}
  \Theta_k(t) 
  = \bigl(v_k+v_{k+1}\bigr) (y_{A_k}-a_k)
  = \cO\bigl(e^{v_kv_{k+1}t}\bigr) \quad(t\to-\infty)
\end{equation}
where $y_{A_k}$ is the $y$-coordinate of $A_k(t)$.  \upshape\medskip

We postpone the elementary proofs of these two Lemmas until \S\ref{sec-proof-of-lem-intersection-location}, \S\ref{sec-proof-of-lem-intersection-angle}.

\subsection{The broken solution}\label{sec-broken-initial-curve}
Assuming the convexity hypothesis \eqref{eq-convexity-hypothesis} we introduce a family of piecewise smooth curves $\Cb_t$ such that $\{\Cb_t\}_{t\le t_0}$ is a solution to Curve Shortening, except at those points where it is not smooth. For simplicity we will describe the construction of the noncompact ancient solution with data $\{a_0,\dots,a_n\}$ and $\{C_1, \dots,C_n\}$. The construction of the compact curve with the same data is very similar.

\begin{figure}[t]
  \includegraphics{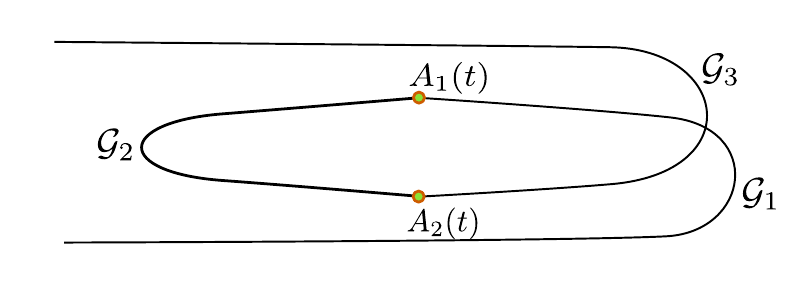}
  \caption{The piecewise smooth convex curve $\Cb(t)$.}
  \label{fig-Cbar}
\end{figure}

We have a collection of translating solitons $\cG_k(t):\;x=G_k(y, t)$, for $k\in\{1, \dots,n\}$ such that $\cG_k(t)$ and $\cG_{k+1}(t)$ share a common asymptote $y=a_k$.  There is a $t_0\in\R$ such that for all $t\le t_0$ the curves $\cG_{k}(t)$ and $\cG_{k+1}(t)$ intersect at some point $A_k(t)$ close to their common asymptote $y=a_{k}$.  We now define the piecewise smooth curve $\Cb_t$ to be the concatenation of the following arcs (see also Figure~\ref{fig-Cbar}):
\begin{itemize}
  \item The part of $\cG_1(t)$ coming in from infinity along the asymptote $y=a_0$ and ending at the intersection point $A_1(t)$,
  \item For $i=2,\dots,n-1$ the segment of $\cG_{i}(t)$ starting at $A_{i-1}(t)$ and ending at $A_{i}(t)$,
  \item the segment of $\cG_{n}(t)$ starting at $A_{n-1}(t)$ and going off to infinity along the asymptote $y=a_n$.
\end{itemize}

\subsection{The really-old-but-not-ancient solutions}
\label{sec-really-old-not-ancient}
The curves $\Cb(t)$ are not smooth so they are not a solution of Curve Shortening.  To construct an actual solution near $\Cb(t)$, we introduce a sequence $\cC_j(t)$ of solutions to Curve Shortening, each defined on some time interval starting at $t=-j$, and with initial value
\[
  \cC_j(-j) = \Cb(-j).
\]
Since $\Cb(-j)$ is a piecewise smooth curve the solution $\cC_j(t)$ exists for $t\in[-j, T_j)$ for some $T_j>-j$.  We will show that
\begin{equation}
  \label{eq-Tj-lower-bound}
  T_* \stackrel{\rm def}= \inf_{j\ge1} T_j > -\infty.
\end{equation}
In addition, we will show that one can extract a subsequence of the solutions $\cC_j(t)$ that converges on any time interval $[T_-, T_*]$ with $T_-<T_*$.  The limit of such a convergent subsequence is then our desired ancient solution.

Our proof of (\ref{eq-Tj-lower-bound}) and the convergence of some subsequence requires two ingredients: the fact that the smooth curves $\cC_j(t)$ ``lie on one side'' of $\Cb(t)$, and an estimate for the area between the two curves.  To explain this in more detail let us first compare the solution $\cC_j(t)$ and the broken solution $\Cb(t)$ for $t>-j$.

\subsection{Orientation of $\Cb(t)$}
\label{sec-orientation-of-Cb}
\begin{figure}[h]
  \centering
  \includegraphics[scale=1.7321]{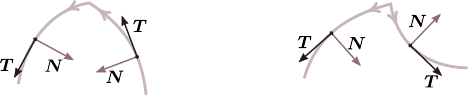}
  \caption{Convex and non-convex corners}
  \label{fig-convex-nonconvex}
\end{figure}

We can orient each smooth arc of the broken solution $\Cb(t)$ by traversing it in such a way that the curve always ``bends to the left.''  More precisely, we choose the unit tangent and normal $\{\vT, \vN\}$ at any smooth location to be a right handed basis for which $\vN$ points in the direction of curvature ($\Cb_{ss} = k\vN$ with $k>0$).  It is clear from the construction in \S~\ref{sec-broken-initial-curve} of the broken curve that all corner points are ``convex'' in the sense that these local orientations match at the corner points (see also figures \ref{fig-intersectionpoint-convex}, \ref{fig-Cbar}, and \ref{fig-convex-nonconvex}).  In any small neighborhood of a point on the curve $\Cb(t)$ we can therefore distinguish unambiguously between points that lie to the left of the curve (i.e.~in the direction of $\vN$) and points to the right of the curve.

\subsection{Separation between $\Cb(t)$ and $\cC_j(t)$}
\begin{figure}[h]
  \centering
  \includegraphics[scale=0.5]{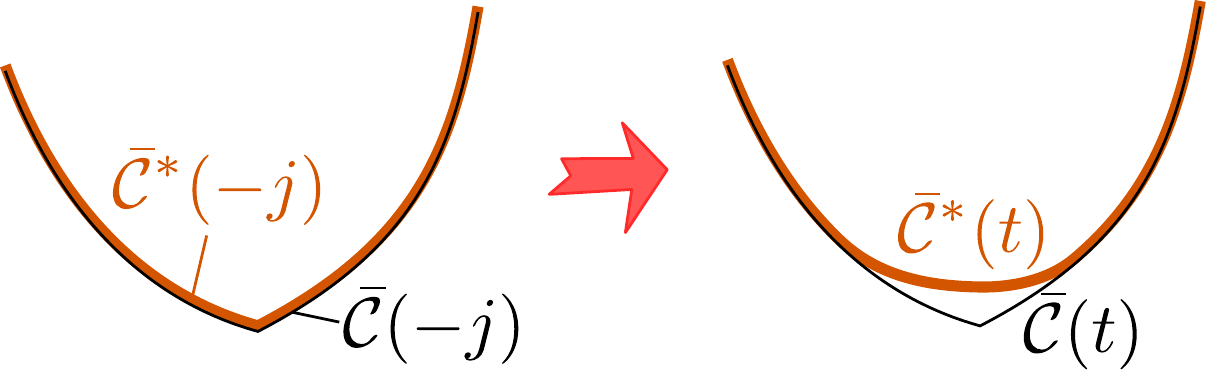}
  \caption{Separation of $\Cb(t)$ and $\cC_j(t)$}
  \label{fig-Cbar-Cj-separation}
\end{figure}

The broken solution $\Cb(t)$ consists of smooth arcs $A_{i-1}(t)A_i(t)$ ($1\le i\le n$) which evolve by Curve Shortening.  Since $\cC_j(t)$ is built out of these arcs at time $t=-j$, it will be very close to $\Cb(t)$ for $t>-j$ close to $-j$.  Initially it is only at the corners $A_i(t)$ where the deviation of $\cC_j(t)$ from $\Cb(t)$ is considerable.  A local analysis of the solutions near any corner point shows that, since the corners are convex, the solution $\cC_j(t)$ will lie to the left of the broken solution $\Cb(t)$.  In fact, the asymptotic shape of the smooth solution $\cC_j(t)$ near any corner $A_i(t)$ is given by a ``Brakke wedge'' (\cite{MR485012,MR3043378}) See Figure~\ref{fig-Cbar-Cj-separation}.

\begin{figure}[h]
  \includegraphics{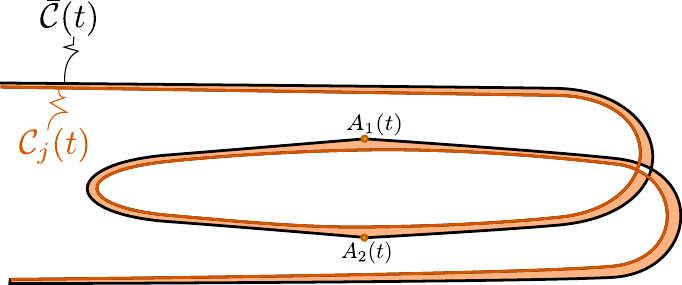}
  \caption{The curves $\Cb(t)$ and $\cC_j(t)$, and the area they enclose.  The area is to be counted with multiplicity.}
  \label{fig-region-between-Cj-Cbar}
\end{figure}

In general, the fact that the corners of $\Cb(t)$ are convex implies that if $\cC^*(t)$ is a solution of Curve Shortening that starts at time $t=-j$ in a narrow strip to the left of $\Cb(-j)$, then $\cC^*(t)$ will remain on the left of $\Cb(t)$ because the maximum principle keeps it from crossing the smooth parts of $\Cb(t)$, and a smooth curve on the left of $\Cb(t)$ cannot touch any of the corner points.  By approximating the initial curve $\cC_j(-j) = \Cb(-j)$ with smooth curves $\cC_{j,\eta}$ that run parallel to $\Cb(-j)$ at a distance no more than $\eta>0$ on the left, and then letting $\eta\to0$, we find that the same holds for $\cC_j(t)$.  Thus it follows that initially, for $t$ shortly after $t=-j$, the whole solution $\cC_j(t)$ will lie in a narrow strip on the left of $\Cb(t)$ (i.e.~in the direction of its normal).  See Figure~\ref{fig-region-between-Cj-Cbar}.

\subsection{The area between $\cC_j(t)$ and $\Cb(t)$}
\label{sec-area-between-Cj-and-Cbar}
For as long as $\cC_j(t)$ is close to $\Cb(t)$, and lies on the left of $\Cb(t)$ we consider the area $A_j(t)$ of the strip enclosed by the two curves.  This area is to be counted with multiplicity as some points may be enclosed more than once.  We define this area to be
\begin{equation}
  \label{eq-area-defined}
  A_j(t) = \int_{\cC_j(t)} xdy - \int_{\Cb(t)} xdy, 
\end{equation}
where both curves are given the orientation described in \S\ref{sec-orientation-of-Cb}.

\subsection{Area growth lemma}
\label{lem-area-growth}
\itshape The area $A_j(t)$ increases according to
\begin{equation}
  \label{eq-area-rate-of-change}
  A_j'(t) = \sum_{k=1}^n \Theta_k(t).
\end{equation}
There exist $T_0<0$ and there exist constants $C, \delta>0$ such that $A_j(t)$ is bounded by
\begin{equation}
  \label{eq-area-bound}
  A_j(t) \le Ce^{\delta t}
\end{equation}
for all $t\le T_0$ and for all $j\ge T_0$.  \upshape \medskip

\textit{Proof. } A family of arcs $\cC: [a, b]\times(t_0, t_1) \to \R^2$ that evolve by Curve Shortening sweeps out area according to
\[
  \frac{dA} {dt} = \int_{\cC}Vds = \int_{\cC}\kappa ds = \bigl[\Delta\theta\bigr]_{\pd\cC}.
\]
In other words, the rate at which the arc sweeps out area is exactly the change in the tangent angle $\theta$ as one goes along the curve.  The rate at which a piecewise smooth solution to Curve Shortening sweeps out area is obtained by summing $[\Delta\theta]_\cC$ over each smooth arc.

For the smooth curves $\cC_j$ the change in $\theta$ is $n\pi$:
\[
  \bigl[\Delta\theta\bigr]_{\cC_j(t)} = n\pi.
\]
Namely, the curve starts at the first asymptote $y=a_0$ and then makes $n$ turns before ending at the last asymptote $y=a_n$.  For the broken curves $\Cb(t)$ the tangent angle also starts at $\theta=0$ and ends at $n\pi$, but along the way $\theta$ makes small jumps of size $\Theta_i(t)$ at each corner point $A_i(t)$.  Thus for the broken solution we get
\[
  \bigl[\Delta\theta\bigr]_{\Cb(t)} = n\pi - \sum_{k}\Theta_k(t).
\]

Subtracting the two we find that the rate at which the area between $\Cb(t)$ and $\cC_j(t)$ grows is as given by \eqref{eq-area-rate-of-change}.

Once we have proved \eqref{eq-area-rate-of-change}, it follows immediately from \eqref{eq-intersection-angle-asymptotics} that $A_j'(t) \le C e^{\delta t}$ for some $C>0$, and for any
\[
  \delta < \min_{1\le k\le n} v_{k-1}v_{k}.
\]Since $A_j(-j) = 0$, we have
\[
  A_j(t) = \int_{-j}^t A_j'(\tau)\, d\tau \le C\int_{-j}^t e^{\delta\tau}\, d\tau \le \frac{C}{\delta} e^{\delta t},
\]
which then leads to the estimate \eqref{eq-area-bound}.

\subsection{Lemma} \itshape
\label{sec-convergence}
For any point $P$ on $\cC_j(t)$ there is a point $Q$ on $\Cb(t)$ such that the distance between $P$ and $Q$ is at most $C e^{\delta t/2}$.  \upshape \medskip

\begin{figure}[h]
  \includegraphics[scale=0.5]{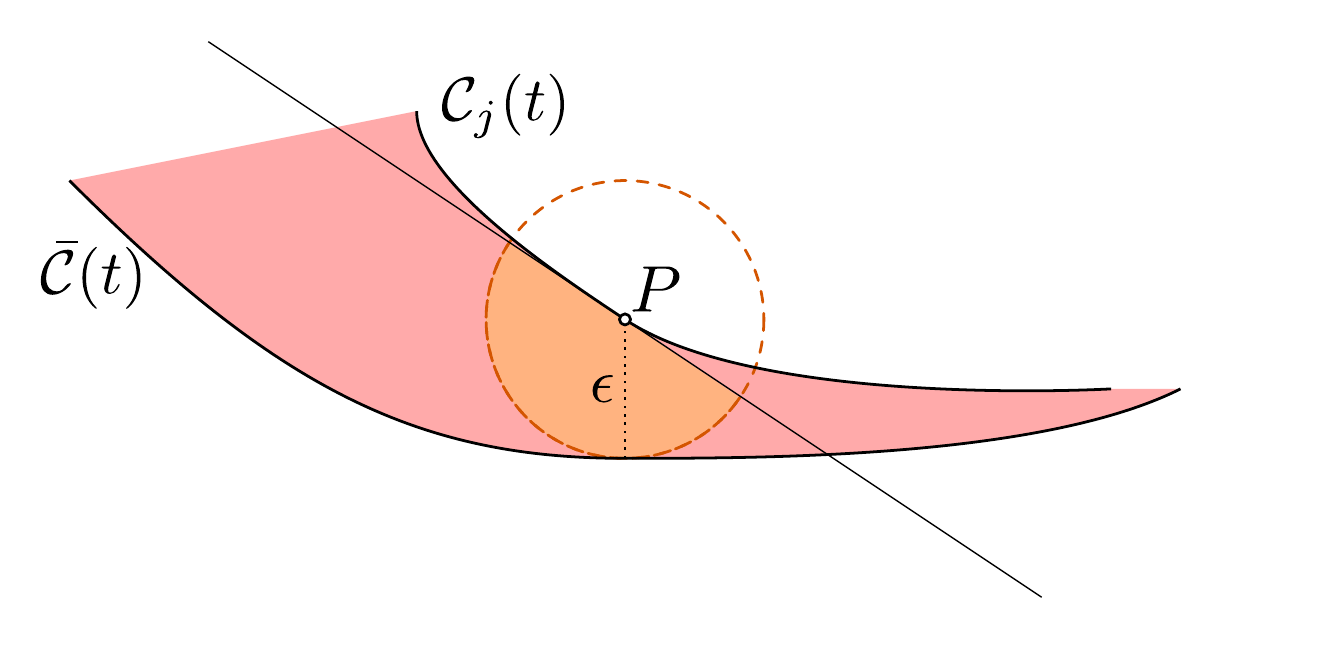}
  \caption{Estimating the distance between two convex curves in terms of the area between them.}
  \label{fig-distance-from-area-estimate}
\end{figure}

This Lemma implies that the curves $\cC_j(t)$ lie within a strip of width $C e^{\delta t/2}$ around $\Cb(t)$.

To prove the Lemma we let $P$ be a given point on $\cC_j(t)$ and let $Q$ be the point on $\Cb(t)$ nearest to $P$.  Then the open disc $B_\epsilon(P)$ is disjoint from $\Cb(t)$, if $\epsilon$ is the distance between $P$ and $Q$.  Moreover, the curve $\cC_j(t)$, being locally convex, must lie on one side of the tangent to $\cC_j(t)$ at $P$.  It follows that the region between $\cC_j(t)$ and $\Cb(t)$ contains at least half of the disc $B_\epsilon(P)$.  See Figure~\ref{fig-distance-from-area-estimate}.  Therefore $A_j(t) \ge \tfrac12 \pi \epsilon^2$.  The Lemma now follows from our estimate \ref{eq-area-bound} for $A_j(t)$.

\subsection{Conclusion of the construction in the convex case}
It is now possible to complete the argument that was started in \S~\ref{sec-really-old-not-ancient}.  Our estimates show that the solution $\cC_j(t)$ lies in a strip of width at most $Ce^{\delta t}$ to the left of the broken solution $\Cb(t)$, at least for as long as $\cC_j(t)$ exists.  The constants $C$ and $\delta$ do not depend on $j$.

We claim that the time $T_j$ at which $\cC_j(t)$ becomes singular is bounded from below by some $T_*$ that does not depend on $j$ (cf.~\eqref{eq-Tj-lower-bound}).  If this were not the case, then according to Grayson's description \cite{grayson:hes86} of nearly singular curves, there would be a sequence of times $t_j\to-\infty$ at which the curve $\cC_j(t_j)$ would have an arc of arbitrarily short length (say, no more than $2^{-j}$) with total curvature arbitrarily close to $\pi$ (say, at least $\pi-2^{-j}$).  This is impossible because $\cC_j(t_j)$ is a smooth convex curve without inflection points which stays in a strip of width $Ce^{\delta t_j}$ to the left of $\Cb(t_j)$.

\section{Construction of embedded ancient non convex solutions}
\label{sec-embedded-ancient-non-convex}
\subsection{Embedded eternal solutions}
\label{sec-eternal-solutions}
We again assume that a finite sequence of heights $a_0, \dots, a_n$ and a corresponding sequence of horizontal shifts $C_1, \dots, C_n$ are given.  Instead of the convexity hypothesis \eqref{eq-convexity-hypothesis} we now assume that the $a_i$ are monotone
\begin{equation}
  \label{eq-monotone-hypothesis}
  a_0 < a_1 < \dots < a_{n-1} < a_n \, .
\end{equation}
These data determine a sequence of Grim Reapers $\cG_k(t)$ as in \eqref{eq-asymptotic-to-GR}.  The monotonicity \eqref{eq-monotone-hypothesis} of the heights $a_i$ implies that the $\cG_k(t)$ are disjoint, and that two consecutive Grim Reapers $\cG_k(t)$ and $\cG_{k+1}(t)$ have the line $y=a_k$ as common asymptote.

We will now construct an ancient solution $\cC(t)$ which lies within the strip $a_0<y<a_n$, and which for $t\to-\infty$ is asymptotic to the Grim Reapers $\cG_k(t)$.  By construction our solution will be a graph $x=u(y, t)$ over the interval $a_0<y<a_n$.  Because of this it will never become singular and is therefore in fact an ``eternal solution.''

The outline of our construction is the same as for the convex ancient solutions from \S~\ref{sec-ancient-convex}: we begin by defining an approximate solution $\Ca(t)$ of Curve Shortening for $t<-T_*$ obtained by gluing together the given family of Grim Reapers.  Then we consider the old-but-not-ancient solutions $\cC_j(t)$ to Curve Shortening that at $t=-j$ start with $\cC_j(t) = \Ca(-j)$.  We then show that the $\cC_j(t)$ are so close to the approximate solutions $\Ca(t)$ that one can extract a convergent subsequence whose limit is an ancient solution.

\begin{figure}[t]
  \centering\def\svgwidth{0.6\textwidth} 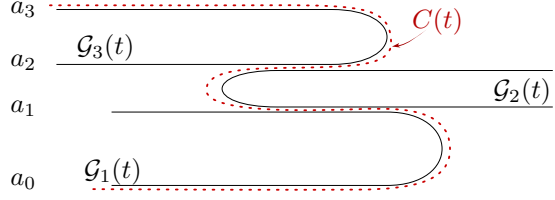
  \caption{An eternal embedded solution $C(t)$.  For $t\ll0$ the curve $C(t)$ is asymptotic to the three Grim Reapers $\cG_i(t)$ ($i=1, 2, 3$.) }
  \label{fig-eternal-embedded}
\end{figure}

\subsection{The approximate solution $\Ca(t)$}
\label{sec-embedded-initial-data}
We can choose $T_*>0$ so large that the Grim Reapers $\cG_k(t)$ ($1\le k \le n$) intersect both vertical lines $x=\pm1$ whenever $t<-T_*$.  For $j > T_*$ we then construct a new curve $\Ca(t)$ by gluing the Grim Reapers along their common asymptote in the strip $-1\le x\le 1$.  To describe the gluing in more detail, we first choose a function $\eta\in C^\infty(\R)$ with $\eta'(x)\ge0$, $\eta(x)=0$ for $x\le -1$, and $\eta(x)=1$ for $x\ge 1$.  Let the Grim Reapers that are asymptotic to $y=a_k$ be given by $y=g_k(x, t)$ and $y=g_{k+1}(x, t)$, respectively, where according to \eqref{eq-grim-reaper},
\begin{equation}
  \label{eq-GR-one-branch}
  \left\{
    \begin{gathered}
      g_{k+1}(x, t) = a_k + \frac{1} {v_{k+1}}\arcsin  e^{\phi_{k+1}}, \quad \phi_{k+1} =  v_{k+1}^2 t + v_{k+1}(x - C_{k+1}),\\
      g_{k}(x, t) = a_k - \frac{1} {v_{k}} \arcsin e^{\phi_k}, \quad \phi_k = v_{k}^2 t - v_{k}(x - C_{k})
    \end{gathered}
  \right.
\end{equation}
(when $k$ is odd; for even $k$ one has to change the sign of $x$.)

\begin{figure}[t]
  \centering \def\svgwidth{0.48\textwidth} 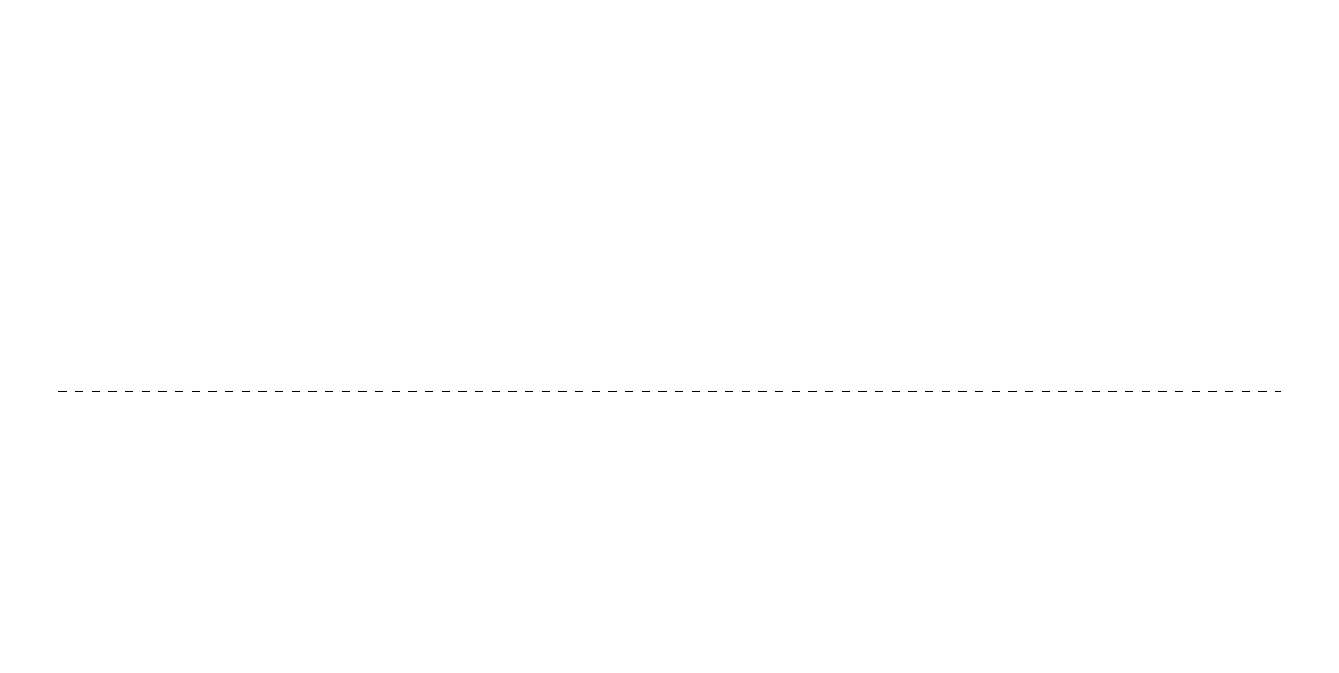 \def\svgwidth{0.48\textwidth} 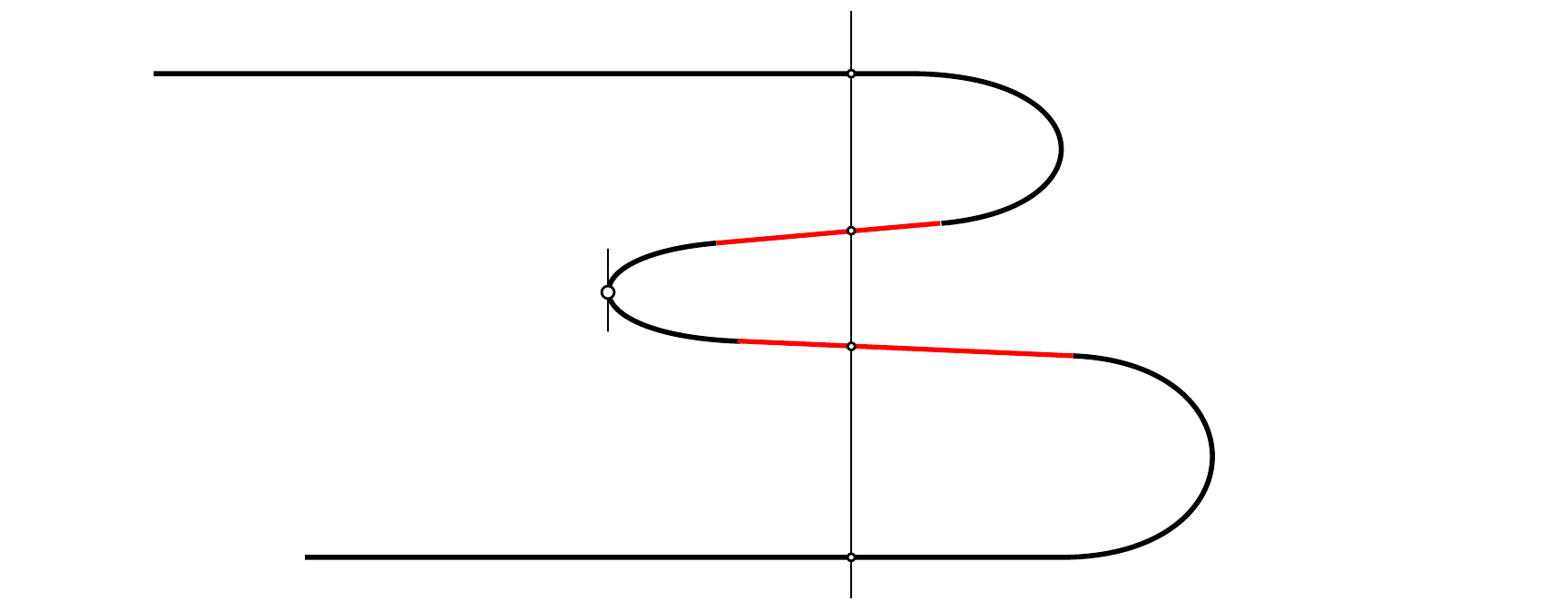
  \caption{\textbf{Left:} gluing two Grim Reapers.  \textbf{Right:} the resulting approximate solution $\Ca(t)$, its intersections with the $y$-axis, and its vertical tangents.}
  \label{fig-embedded-initial-curve}
\end{figure}
We can then replace the two Grim Reaper branches $y=g_{k}(x, t)$, $g_{k+1}(x, t)$ by
\[
  \bar g(x,t) = \eta(x) g_{k+1}(x, t) + \bigl(1-\eta(x)\bigr) g_k(x, t).
\]
It is easy to verify that
\[
  g_k(x, t) \le \bar g(x, t) \le g_{k+1}(x, t), \text{ and } \bar g_x(x, t)>0.
\]
We repeat this procedure for each $k=1, \dots, n-1$ which, in the end, results in a curve $\Ca(t)$ with the following properties for all $t<-T_*$:
\begin{enumerate}
  \item $\Ca(t)$ intersects the $y$-axis in exactly $n+1$ points $P_0$, $P_1$, \dots, $P_n$
  \item in between each pair of consecutive intersection points $P_{k-1}$, $P_{k}$ the curve $\Ca(t)$ has exactly one point $Q_k$ with a vertical tangent
  \item the segments of the curve between $Q_k$, $Q_{k+1}$ are graphs of functions $y=u_k(x, t)$.
\end{enumerate}

\subsection{The old-but-not-ancient solutions $\cC_j(t)$}
\label{sec-embedded-old-but-not-ancient}
For any $j<-T_*$ we let $\{\cC_j(t) : t\ge -j\}$ be the solution of Curve Shortening which at time $t=-j$ starts with $\cC_j(-j) = \Ca(-j)$.  Since $\Ca(-j)$ is the graph of a function $x=f(y)$ defined on the interval $a_0<y<a_n$, it is easy to show that the solution $\cC_j(t)$ never becomes singular and exists for all $t>-j$.

For $t>-j$ the maximum principle implies that the solution $\cC_j(t)$ must still avoid the Grim Reapers $\cG_k(t)$ ($k=1, \dots, n$).  Therefore $\cC_j(t)$ intersects the $y$-axis in $n$ distinct points $P_{jk}(t)$ ($k=1,\dots,n$), and in between each pair $P_{j, k-1}(t), P_{j, k}(t)$ there is a unique point $Q_{j, k}(t)$ with a vertical tangent.  Each arc $Q_{j, k}(t)Q_{j, k+1}(t)$ is the graph of a function $y=u_{jk}(x, t)$.

\begin{figure}[h]
  \centering
  \includegraphics[width=0.66\textwidth]{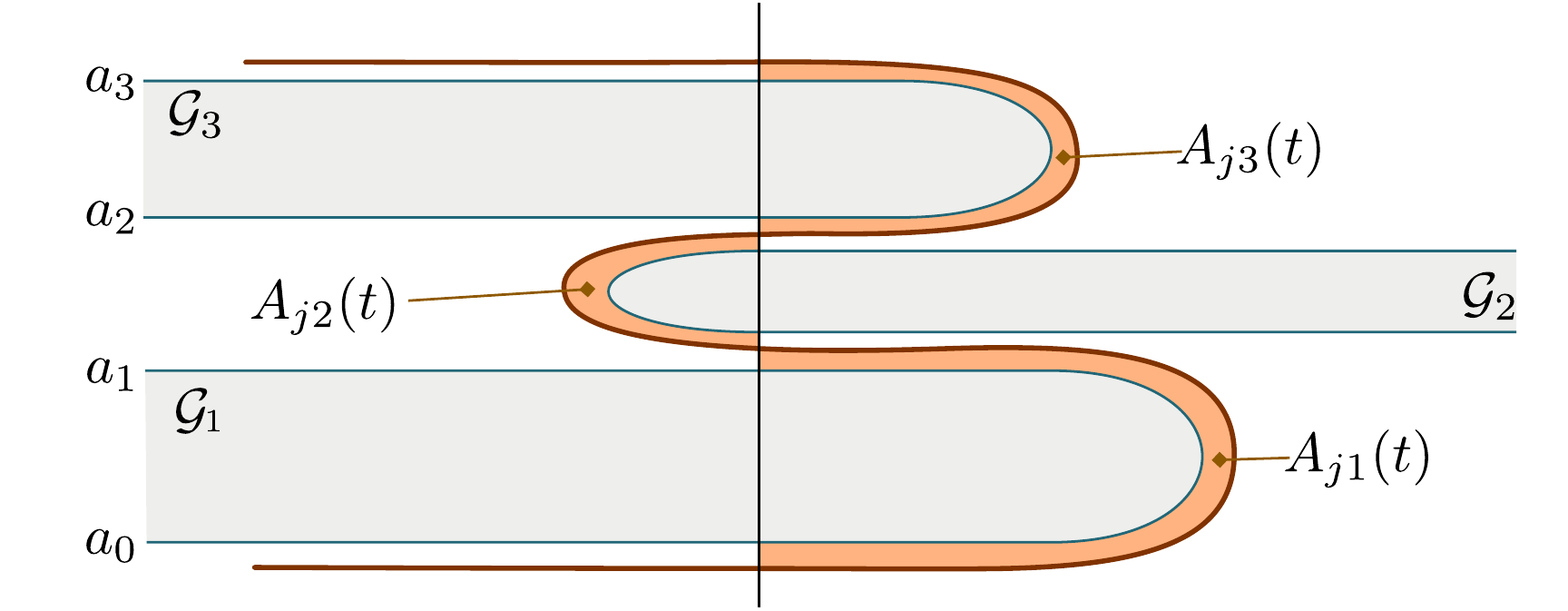}
  \hfill
  \includegraphics[width=0.3\textwidth]{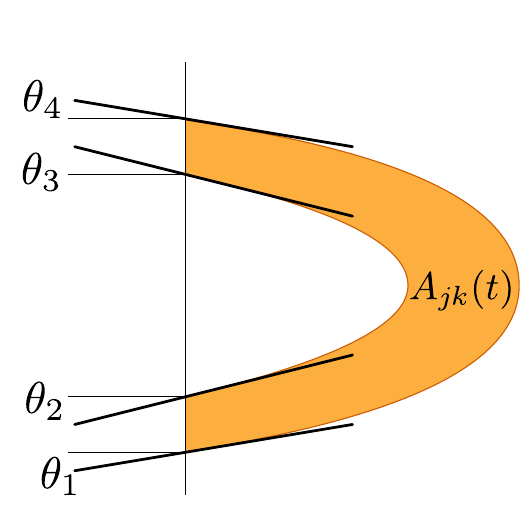}
  \caption{\textbf{Left: }The areas $A_{jk}(t)$ between the Grim Reapers $\cG_{k}(t)$ and the solutions $\Cb_j(t)$.  \textbf{Right:} The angles $\theta_i$ corresponding to one of the areas $A_{jk}(t)$.}
  \label{fig-areas}
\end{figure}

\subsection{The areas $A_{jk}(t)$}
As in the convex case we show that the old-but-not-ancient solutions $\cC_j(t)$ are close to the Grim Reapers $\cG_k(t)$ by estimating the area between them.  To be precise, we let $A_{jk}(t)$ be the area of the horseshoe shaped region $A_{jk}(t)$ enclosed by the $y$-axis, the Grim Reaper $\cG_k(t)$, and the arc $P_{j,k-1}(t)P_{j, k}(t)$ of $\Cb(t)$ (see Figure~\ref{fig-areas}).  The four edges of the region whose area is $A_{jk}(t)$ are
\begin{itemize}
  \item a segment of the Grim Reaper $\cG_k(t)$,
  \item the segment $P_{j, k-1}(t)P_{jk}(t)$ of the old-but-not-ancient solution $\cC_j(t)$, and
  \item two short segments on the $y$-axis.
\end{itemize}
Thus the edges of $A_{jk}(t)$ all move by Curve Shortening so that we know that
\begin{equation}
  \label{eq-area-rate-of-change-2}
  \frac{d A_{jk}(t)} {dt} = \int_{\pd A_{jk}} \kappa ds
  = -\theta_1 + \theta_2 + \theta_3 - \theta_4
\end{equation}
where $\theta_i$ are the angles indicated in Figure~\ref{fig-areas}.

\subsection{Lemma}\itshape%
There is a constant $M<\infty$ depending only on the $a_i$ and $C_j$ such that
\begin{equation}
  |\theta_1|+|\theta_2|+|\theta_3|+|\theta_4| \le Me^{\delta t}
  \label{eq-angle-estimate}
\end{equation}
for all $t\in[-j, T_*]$.

\upshape\medskip

\begin{proof}
The angles $\theta_2$ and $\theta_3$ are the angle of intersection of the Grim Reaper $\cG_k(t)$ with the $y$-axis.  The branches of the two Grim Reapers that are asymptotic to the line $y=a_k$ are given by \eqref{eq-GR-one-branch}, from which one directly finds an exponential upper bound the form $Me^{\delta t}$ for the slope at $x=0$, and thus for the angles $\theta_2$ and $\theta_3$.
  
The solution $\cC_j(t)$ near its intersection point $P_{jk}(t)$ is a graph $y=u_{jk}(x,t)$ where $u$ satisfies
\begin{equation}
  \label{eq-CS-graphs}
  u_t = \frac{u_{xx}} {1+u_x^2}.
\end{equation}
Since the graph of $u_{jk}$ is caught between the two Grim Reapers $\cG_k(t)$ and $\cG_{k+1}(t)$, the explicit representation \eqref{eq-GR-one-branch} implies that $u$ is bounded by
\[
  |u(x, t) - a_k| \le N e^{\delta t}\text{ for }|x|\le 1, -j\le t <T_*
\]
for some $\delta>0$ and $N<\infty$.  The interior gradient estimates for graphical Mean Curvature Flow (see Evans-Spruck \cite[\S5.2]{MR1151756}) now imply that $|u_x|$ is uniformly bounded on a smaller interval, i.e.~there is some constant $N'<\infty$ such that
\[
  |u_x(x,t)| \le N' \text{ for } -\tfrac12 < x < \tfrac12, -j\le t\le T_*.
\]%
We then know that the parabolic equation \eqref{eq-CS-graphs} for $u$ is nondegenerate on $(-\frac12, \frac12)\times (-j, T_*)$, so that standard parabolic estimates \cite[chapter VI]{ladyzhenskaia1988linear} imply
\[
  |u_x(0,t)| \le N'' e^{\delta'' t}, \text{ for }-j\le t \le T_*
\]
again, for certain constants $N''<\infty$ and $\delta''>0$.
\end{proof}

\subsection{Lemma}\itshape\label{lem-embedded-area-estimate}%
The areas $A_{jk}(t)$ are uniformly bounded by
\[
  A_{jk}(t) \le M e^{\delta t} \text{ for } -j\le t \le T_*
\]
for all $j$ and for $k=1, \dots, n$.  \upshape\medskip

\begin{proof}
The construction of the initial curve $\cC_j(-j)$ is such that at time $t=-j$ the area $A_{jk}(-j)$ is contained in the region $|x|\le 1$, $g_{k}(x, -j)\le y\le g_{k+1}(x, -j)$.  On the interval $-1\le x\le 1$ the explicit expressions \eqref{eq-GR-one-branch} imply that
\[
  \left|g_{k+1}(x, -j) - g_k(x, -j)\right| \le Me^{-\delta j} = M e^{\delta t}
\]
for suitably chosen $M<\infty$ and $\delta>0$.  Hence the areas $A_{jk}(-j)$ satisfy a similar estimate.
  
The Lemma follows by integrating \eqref{eq-area-rate-of-change-2} and using the upper bounds \eqref{eq-angle-estimate}.
  
\end{proof}

\begin{figure}[h]\centering
  \includegraphics[width=0.5\textwidth]{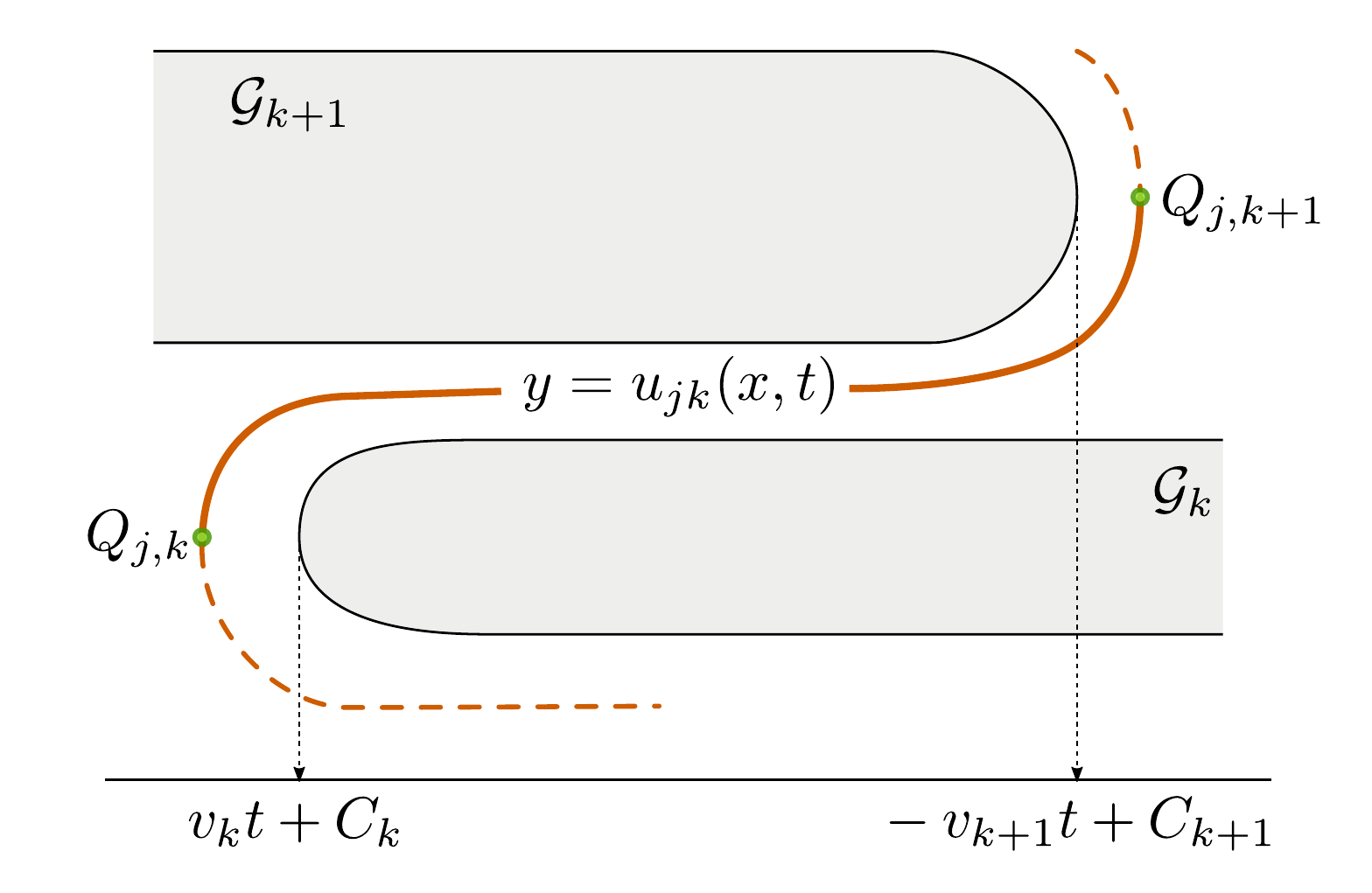}
\end{figure}

\subsection{The arcs $Q_{j,k}Q_{j,k+1}$}

Consider any one of the arcs $Q_{j,k}(t)Q_{j, k+1}(t)$ on $\cC_j(t)$ between to consecutive vertical tangents.  This arc is a graph $y=u_{jk}(x, t)$ (see \S~\ref{sec-embedded-old-but-not-ancient}), where $u_{jk}$ satisfies the Curve Shortening for graphs equation \eqref{eq-CS-graphs}.  To simplify notation we will write $u$ instead of $u_{jk}$ in this section.  Since the arc lies between $\cG_k(t)$ and $\cG_{k+1}(t)$, the function $u(x, t)$ is defined on an interval containing the $x$-coordinates of the tips of the two Grim Reapers, i.e.~$u(x,t)$ is defined for
\[
  -\infty <t\le -T_*, \quad v_kt + C_k < x < -v_{k+1}t + C_{k+1},
\]
and on this domain $u$ is bounded by~\eqref{eq-GR-one-branch}, i.e.
\[
  g_k(x, t) < u(x, t) < g_{k+1}(x, t).
\]
Furthermore, the area bounds from Lemma~\ref{lem-embedded-area-estimate} imply
\begin{equation}
  \label{eq-embedded-L1-bound-a}
  \int_{v_kt+C_k}^0 \bigl( u(x, t) - g_k(x,t) \bigr) dx\le M e^{\delta t}
\end{equation}
and
\begin{equation}
  \label{eq-embedded-L1-bound-b}
  \int_0^{-v_{k+1}t+C_{k+1}} \bigl( g_{k+1}(x, t) - u(x, t)\bigr) dx \le M e^{\delta t}
\end{equation}
for all $t\in[-j, -T_*]$ (note that the integrands are positive).

\subsection{Lemma} \label{sec-embedded-ujk-estimate}\itshape%
For any $L>1$ there are constants $M, \delta>0$ such that for $-j\le t\le -T_*$ one has
\[
  \bigl|\pd_x^mu(x, t) - \pd_x^m g_k(x, t)\bigr| \le Me^{\delta t} \text{ for }1 \le x-v_kt-C_k \le L
\]
and
\[
  \bigl|\pd_x^mu(x, t) - \pd_x^m g_{k+1}(x, t)\bigr| \le Me^{\delta t} \text{ for } -L \le x+v_{k+1}t-C_{k+1} \le -1
\]
for any $m=0,1,2,\dots$ \upshape\medskip

\begin{proof}
We consider the first case and estimate $u-g_k$.  Both $u$ and $g_k$ are solutions of \eqref{eq-CS-graphs}.  The Evans-Spruck interior estimates for graphical Mean Curvature Flow \cite{MR1151756} imply that $u_x$ is uniformly bounded for $1 \le x-v_kt-C_k \le L$.  Interior regularity for general quasilinear parabolic equations then implies that all higher derivatives $\pd_x^mu$ are also bounded.  Since we have shown that $u-g_k$ is small in $L^1$ (see \eqref{eq-embedded-L1-bound-a} and \eqref{eq-embedded-L1-bound-b}), an interpolation inequality leads to the estimates in the Lemma.
\end{proof}

\begin{figure}[h]
  \def\svgwidth{0.6\textwidth} \input eternalEmbeddedInterceptsVerticalPoints.pdf_tex
\end{figure}

\subsection{Convexity of the tips}
\label{sec-embedded-convex-tips}\itshape%
For each $k$ the section of the arc $P_{k-1}P_k$ (see Figure~\ref{fig-embedded-initial-curve}) on which $x<v_kt+C_k+L$ for even $k$ (or $x>-v_{k}t+C_k-L$ for $k$ odd) is convex.  \upshape\medskip

\begin{proof}
Assume for simplicity that $k$ is even.  Then we have shown in Lemma~\ref{sec-embedded-ujk-estimate} that $u(x, t)$ is $C^2$ close to the Grim Reaper $\cG_k$, at least in the strip $1\le x - v_kt -C_k \le L$.  Since the Grim Reaper is convex, the part of the arc $P_{k-1}P_{k}$ on which $1\le x - v_kt -C_k \le L$ must also be convex.  To verify that the remaining part, where $x-v_kt-C_k\le1$, is also convex, we note that this is certainly so at time $t=-j$ because $\cC_j(-j) = \Ca(-j)$ and $\Ca(-j)$ coincides with a Grim Reaper for $x-v_kt-C_k\le1$.  For $t>-j$ we recall that the curvature $\kappa_j$ satisfies a parabolic equation $\kappa_t = \kappa_{ss}+\kappa^3$, so that the maximum principle forces $\kappa>0$ in the region $x-v_kt-C_k\le1$.
\end{proof}

\subsection{Lemma} \label{lem-embedded-uniformly-close}\itshape%
There exist $N, \delta>0$ so that each $\cC_j(t)$ lies within a strip of width $Ne^{\delta t}$ around $\Ca(t)$.  The curvature of $\cC_j(t)$ is uniformly bounded for all $t\in [-j, -T_*]$.  \upshape\medskip

\begin{proof}
Away from the tips, where $x=\pm v_kt+C_k$, we have already established the exponential closeness to the Grim Reapers, and hence to the approximate solution in Lemma~\ref{sec-embedded-ujk-estimate}.  Near the tips, where $0\le x-v_kt-C_k\le L$ (for even $k$), we know that the solution $\cC_j(t)$ and the Grim Reaper $\cG_k(t)$ both are convex, while the area between them is exponentially small.  The same argument as in Lemma~\ref{sec-convergence} then implies that the part of $\cC_j(t)$ with $a_{k-1}<y<a_k$ and $x-v_kt-C_k\le L$ must lie in an $Ne^{\delta t}$ neighborhood of $\cG_k(t)$.
  
The curvature bounds can also be proved by considering the tips and the arcs between them separately.
  
The arcs $Q_kQ_{k+1}$ are graphs of solutions $y=u(x, t)$ to \eqref{eq-CS-graphs}.  Lemma~\ref{sec-embedded-ujk-estimate} provides uniform curvature bounds whenever $v_kt+C_k +1 \le x \le -v_{k+1}t+C_{k+1}$ (for $k$ even).  Any point on a tip region, where $x<v_kt+C_k+L$, lies on a graph $x=v(y, t)$ of a solution to
\[
  v_t = \frac{v_{yy}}{1+v_y^2}
\]
which is bounded by
\[
  v_kt+C_k - Ne^{\delta t} \le v(y, t) \le v_kt+C_k +L.
\]
Interior estimates then again imply that all derivatives $\pd_y^mv$ are uniformly bounded, and in particular, that the curvature $\kappa$ is bounded in the tip region.
\end{proof}

\section{Construction of general ancient non convex solutions}
\label{sec-general-ancient-non-convex}
So far we have shown how one can glue together Grim Reapers so as to construct ancient solutions that are either convex or that are embedded.  In this sections we show how one can combine the two constructions to produce ancient solutions by gluing any finite chain of Grim Reapers with matching asymptotes.

\subsection{Hypotheses in the general case}
\label{sec-general-hypotheses}
\begin{figure}[h]
  \includegraphics[width=0.9\textwidth]{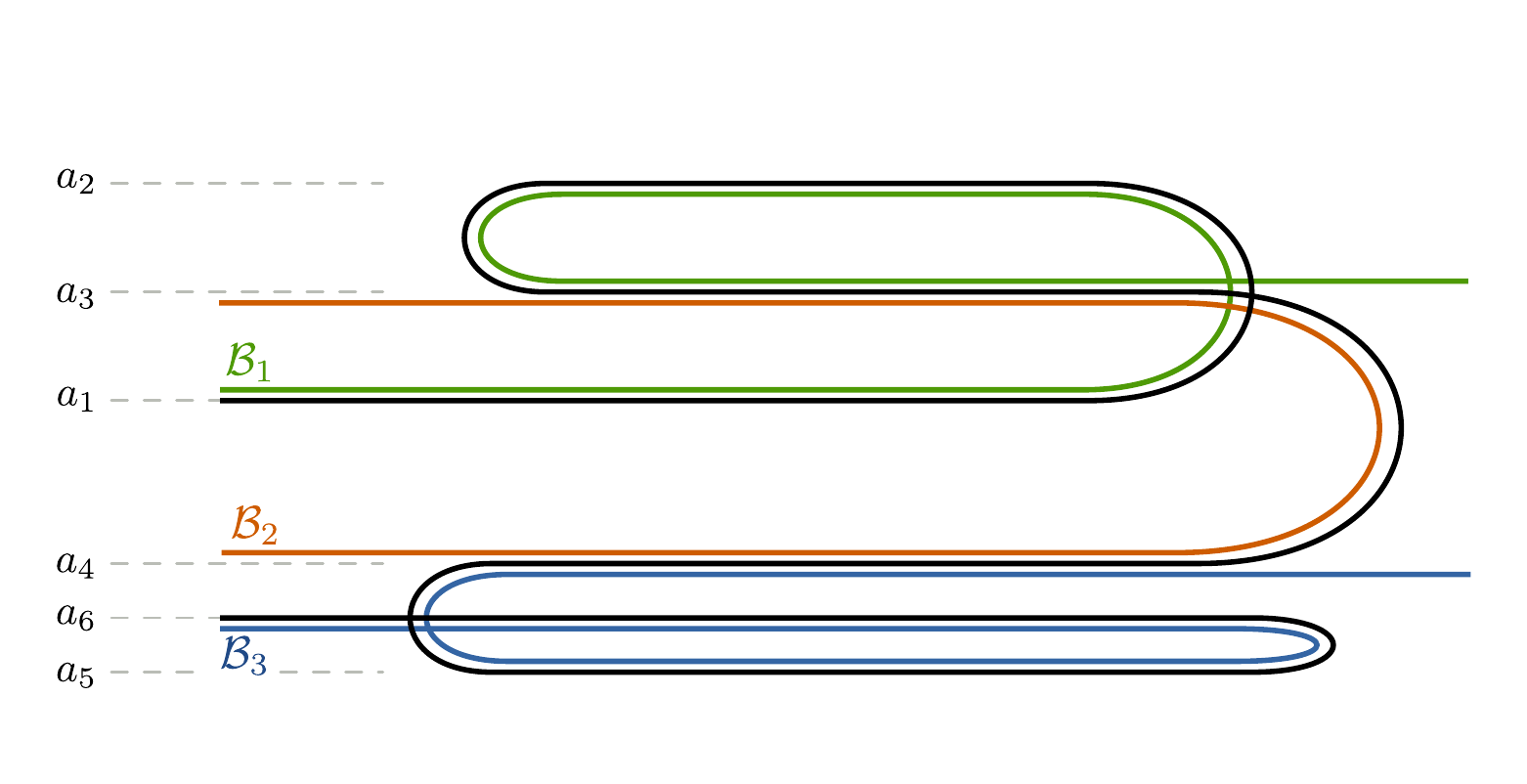}
  \caption{The barriers $\cB_i$ corresponding to a general ancient non convex solution.}
  \label{fig-general-case-barriers}
\end{figure}

Let $\{a_0, a_1, \dots, a_n\}$ be a given collection of heights, and let $\{C_1, \dots, C_n\}$ be a corresponding set of horizontal displacements.  While we are not assuming that the heights $a_j$ alternate as in \eqref{eq-convexity-hypothesis}, we split the sequence $a_j$ into a number of maximal alternating runs.  Thus there are
\[
  0=\ell_0<\ell_1<\ell_2<\dots<\ell_m=n
\]
such that each of the subsequences
\[
  \{a_{\ell_i}, a_{\ell_i + 1}, \dots , a_{\ell_{i+1}-1}, a_{\ell_{i+1}}\}
\]
is a maximal subsequence of $\{a_0, \dots, a_n\}$ that does satisfy (\ref{eq-convexity-hypothesis}) (See Figure~\ref{fig-general-case-barriers}.)

The cases we have considered in the sections~\ref{sec-ancient-convex} and ~\ref{sec-embedded-ancient-non-convex} are in a sense the extreme cases.  On one hand, in the convex case the whole sequence $\{a_0, \dots, a_n\}$ satisfies \eqref{eq-convexity-hypothesis} so that there is only one maximal alternating subsequence; on the other hand, in the embedded case no subsequence of length more than two satisfies \eqref{eq-convexity-hypothesis}, so that every consecutive pair $\{a_j, a_{j+1}\}$ is a maximal alternating subsequence.

Each of the maximal alternating subsequences $\{a_{\ell_i}, \dots, a_{\ell_{i+1}}\}$ satisfies (\ref{eq-convexity-hypothesis}), so there is an ancient convex solution $\cB_i(t)$ that is asymptotically described by the concatenation of the Grim Reapers $\cG_{\ell_i}(t)$, \dots, $\cG_{\ell_{i+1}}(t)$, where $\cG_j(t)$ is given as before in \eqref{eq-grim-reaper}, with $C=C_j$, $a=a_j$, $v=v_j = \pi/|a_{j+1}-a_j|$.  To find an ancient solution that is asymptotic to all Grim Reapers $\cG_1(t)$, \dots, $\cG_n(t)$ we modify the construction of the embedded ancient solitons from section \ref{sec-embedded-ancient-non-convex} by replacing the Grim Reapers $\cG_j(t)$ by the convex ancient solutions $\cB_i(t)$.  As before we can use a cut-off function to glue the barriers $\cB_i(t)$ together and thus create an approximate solution $\Ca(t)$ (see section~\ref{sec-embedded-initial-data}).

We define a sequence of old-but-not-ancient solutions $\cC_j(t)$ ($t\in[-j, T_j)$) by requiring the initial value $\cC_j(-j)$ of $\cC_j$ to come from the approximate solution $\Ca(-j)$.

\begin{figure}[h]
  \includegraphics[width=0.5\textwidth]{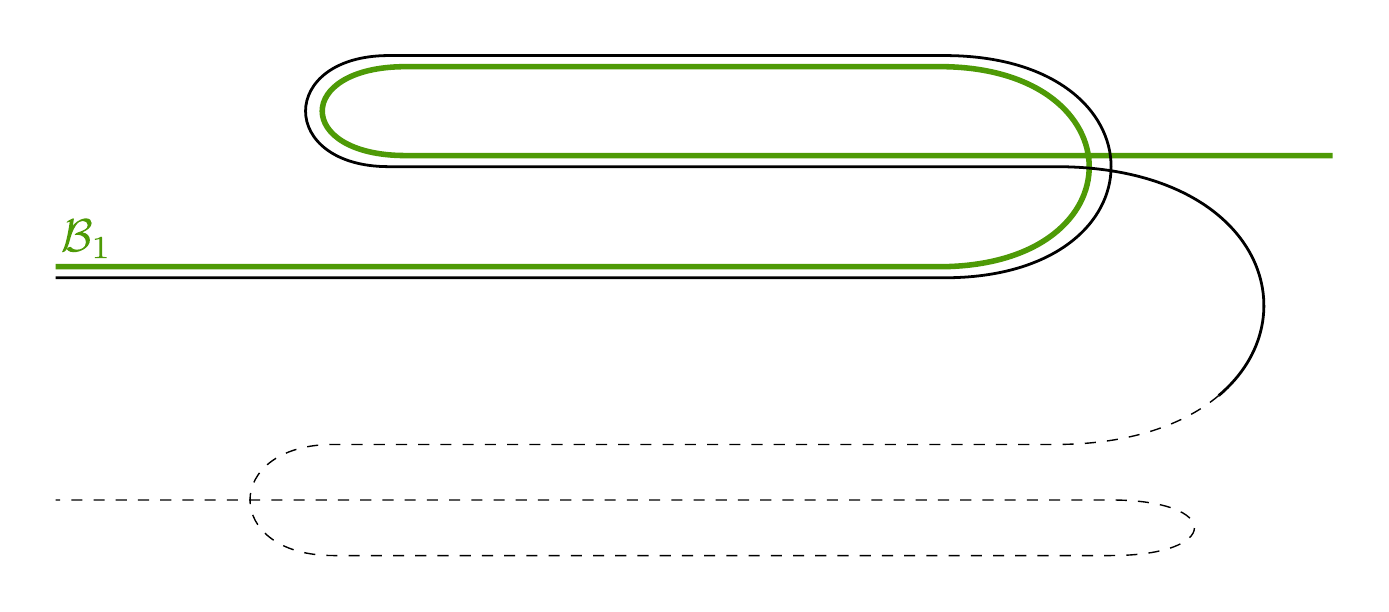}
  \caption{The solutions $\cC_j(t)$ lie on one side of the barriers $\cB_i(t)$.}
  \label{fig-general-case-barriers-one-side}
\end{figure}

The barriers $\cB_i(t)$ are convex, so the solutions $\cC_j(t)$ will locally lie on one side of the $\cB_i(t)$.

The situation in the embedded case, as described in Figure~\ref{fig-embedded-initial-curve} applies here too.  For all $t<-T_*$ the approximate solution $\Ca(t)$ intersects the $y$-axis exactly $n+1$ times (once for each asymptote $y=a_j$).  The solutions $\cC_j(t)$ will initially also have such intersections.  Since they must lie on one side of the barriers, they cannot lose those intersections, and thus we know that the $\cC_j(t)$ intersect the $y$-axis at $n$ distinct points $P_{j1}(t)$, \dots, $P_{jn}(t)$.  In between any consecutive pair $\{P_{j, k}(t), P_{j, k+1}(t)\}$ of intersections there is again exactly one point $Q_{jk}(t)$ with a vertical tangent; the arcs $Q_{j,k-1}(t)Q_{j,k}(t)$ connecting the vertical tangents are graphs of functions $y=v_{j,k}(x, t)$.

\begin{figure}[h]
  \centering
  \includegraphics[width=0.7\textwidth]{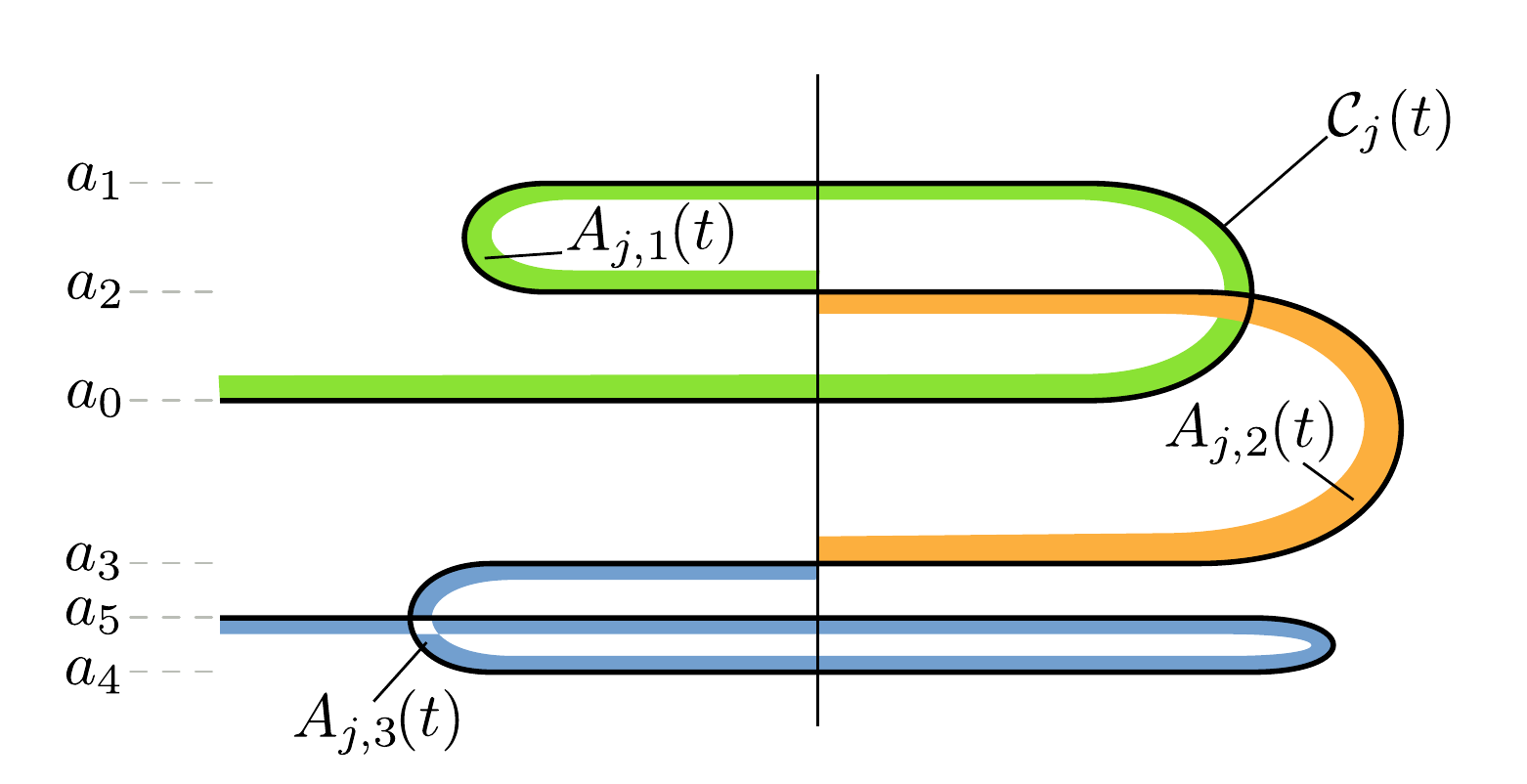}
  \caption{The areas we estimate in the construction of an ancient solution with data $\{a_0, a_1, \dots, a_n; C_1, \dots, C_n\}$.
    Each area is bounded by a segment of the old-but-not-ancient solution $\cC_j(t)$ (black), and a segment of the corresponding convex barrier $\cB_i(t)$.}
  \label{fig-general-case-areas}
\end{figure}

\subsection{The areas $A_{j,k}(t)$}
\label{sec-general-case-the-areas}
The essential point in the construction of the embedded ancient solutions was the area estimate in Lemma~\ref{lem-embedded-area-estimate} which allowed us to conclude that the old-but-not-ancient solutions $\cC_j(t)$ were close to the approximate solutions $\Ca(t)$.  Here we can use a similar argument provided we adapt the definition of the regions whose area we bound.

For each maximal alternating subsequence $\{a_{\ell_i}, \dots, a_{\ell_{i+1}}\}$ we define a region $\cA_{ji}(t)$ which is bounded by the arc $P_{j,\ell_i}(t)P_{j,\ell_{i+1}}(t)$ of $\cC_j(t)$, the corresponding segment of the barrier $\cB_{i}(t)$, and the two short segments on the $y$-axis that connect these two arcs.  See figure~\ref{fig-general-case-areas}.  The growth of the areas of the $\cA_{ji}(t)$ is again determined by the angles at which $\cC_j(t)$ and $\cB_i(t)$ meet the $y$-axis.  Near those intersections both $\cB_i(t)$ and $\cC_j(t)$ are close to either the line $y=a_{\ell_i}$ or the line $y=a_{\ell_{i+1}}$, with exponential bounds like those in Lemma~\ref{sec-embedded-ujk-estimate}.  It follows that the areas $A_{ji}(t)$ area uniformly bounded by $Ne^{\delta t}$ for some $N<\infty$ and $\delta>0$.

\subsection{Curvature bounds and existence of the ancient solution}
To complete the existence proof of an ancient solution we must show that the old-but-not-ancient solutions $\cC_j(t)$ are defined for all $t\in[-j, -T_*]$, for some $T_*$ that does not depend on $j$, and that the curvature of $\cC_j(t)$ is uniformly bounded for all $j$ and all $t\le T_*$.  To this end we can merely repeat the arguments in the proof of Lemma~\ref{lem-embedded-uniformly-close}, so that we are done.

\section{Miscellaneous proofs}
\subsection{Proof of Lemma~\ref{lem-intersection-location}}
\label{sec-proof-of-lem-intersection-location}
Without loss of generality we may assume that $a_{k\pm 1} > a_k$, and that
\begin{gather*}
  G_k(y, t) = v_kt+C_k+\frac{1}{v_k}G\bigl(v_k(y-a_k)\bigr), \\
  G_{k+1}(y, t) = -v_{k+1}t-C_{k+1}-\frac{1}{v_{k+1}}G\bigl(v_{k+1}(y-a_k)\bigr).
\end{gather*}
Here we have used $G(s) = G(\pi-s)$ to arrive at the expression for $G_{k}$.

Let $b_k = a_k + \frac12\min\{|a_{k-1}-a_k|, |a_{k+1}-a_k|\}$.

Existence of the intersection between the two Grim Reapers $x=G_k(y, t)$ and $x=G_{k+1}(y, t)$ in the interval $(a_k, b_k)$ follows from the fact that as $y\searrow a_k$ one has $G_k(y, t)\to+\infty$ and $G_{k+1}(y, t) \to -\infty$, while for $t\to-\infty$ one also has
\begin{gather*}
  G_k(b_k, t) = C_k+\frac{1}{v_k}G(v_k(b_k-a_k)) +v_kt \to -\infty,\\
  G_{k+1}(b_k, t) = -C_{k+1}-\frac{1}{v_{k+1}}G(v_{k+1}(b_k-a_k)) - v_{k+1}t \to +\infty.
\end{gather*}
Thus $G_{k+1}(y, t) - G_k(y, t)$ changes sign on the interval $(a_k, b_k)$.

Uniqueness of the intersection follows from the fact that $G_k$ is decreasing while $G_{k+1}$ is increasing on the interval $(a_k, b_k)$ .

To find the asymptotic location of the intersection we solve the two equations $x=G_{k}(y,t)$, $x=G_{k+1}(y,t)$, resulting in
\begin{subequations}
  \begin{gather}
    \label{eq-intersection-derivation-1}
    y - a_k = \frac{1}{v_k}\arcsin e^{v_k (x-v_kt-C_k)}\\
    y - a_k = \frac{1}{v_{k+1}}\arcsin e^{-v_{k+1} (x+v_{k+1}t+C_{k+1})}.
  \end{gather}
\end{subequations}
As $t\to-\infty$, the two Grim Reapers separate and their intersection approaches their common asymptote $y=a_k$, so we may assume that $\eta \stackrel{\rm def}= y-a_k\to0$.  We then get
\begin{gather*}
  \eta + \cO(\eta^3) = \frac{1}{v_k} e^{v_k(x-v_kt-C_k)},\\
  \eta + \cO(\eta^3) = \frac{1}{v_{k+1}} e^{-v_{k+1}(x+v_{k+1}t+C_{k+1})}.
\end{gather*}
Use $\eta+\cO(\eta^3)= \eta(1+\cO(\eta^2)) = \eta e^{\cO(\eta^2)}$,
\[
  v_k(x-v_kt-C_k) -\ln v_k = -v_{k+1}(x+v_{k+1}t+C_{k+1}) - \ln v_{k+1} + \cO(\eta^2).
\]
and rearrange terms
\[
  (v_k+v_{k+1}) x = \bigl(v_k^2 - v_{k+1}^2\bigr)t + v_kC_k - v_{k+1}C_{k+1} +\ln \frac{v_k}{v_{k+1}} + \cO(\eta^2).
\]
Finally we find
\begin{equation}
  \label{eq-intersection-derivation-2}
  x = \bigl(v_k-v_{k+1}\bigr) t + c_k + \cO(\eta^2)
\end{equation}
with
\[
  c_k = \frac{v_kC_k - v_{k+1}C_{k+1} + \ln v_k/v_{k+1}}{v_k+v_{k+1}}.
\]
Substitute $x-v_kt = -v_{k+1}t + c_k + \cO(\eta^2)$ in \eqref{eq-intersection-derivation-1} to get
\[
  \eta+\cO(\eta^3) = \frac{1}{v_k} e^{v_kv_{k+1}t+c_k'} \bigl(1+\cO(\eta^2)\bigr), \text{ i.e. } \eta = \cO\bigl(e^{v_kv_{k+1}t}\bigr),
\]
for some constant $c_k'$.  Together with \eqref{eq-intersection-derivation-2} this implies \eqref{eq-intersection-convex-x-coord}.

\subsection{Proof of Lemma~\ref{lem-intersection-angle}}
\label{sec-proof-of-lem-intersection-angle}
The translating soliton $y=G_k(x, t)$ has equation
\[
  \sin (v_k\eta) = e^{v_k(x-v_kt -C_k)},
\]
where $\eta = y-a_k = G_k(x, t)-a_k$.  Differentiating w.r.t.~$x$ we find
\[
  v_k\cos (v_k\eta) \frac{\dd y} {\dd x} = v_k e^{v_k(x-v_kt -C_k)} =v_k \sin (v_k\eta) \implies \frac{\pd G_k} {\pd x}(x, t) = \tan(v_k\eta).
\]
It follows that the angle between the tangent to $y=G_k$ and the $x$-axis is exactly $v_k\eta$.  A similar calculation shows that the angle between the $x$-axis and the tangent to $y=G_{k+1}(x, t)$ at $A_k(t)$ is exactly $v_{k+1}\eta$.  The angle the two translating solitons make at their intersection point $A_k(t)$ is therefore $\Theta_k(t) = (v_k+v_{k+1})\eta$.  This implies \eqref{eq-intersection-angle-asymptotics}.

\subsection{Proof of Lemma~\ref{lem-area-growth}}
\label{sec-proof-of-area-growth}
We defined the area between the two locally convex curves $\cC_j(t)$ and $\Cb(t)$ in \eqref{eq-area-defined} by means of a line integral.  The line integral is taken over a piecewise smooth curve (the smooth arcs comprising $\Cb(t)$ and $\cC_j(t)$).  If $\cA_i(t)$ is one of those arcs, and if we parametrize it by $(u,t)\mapsto \bigl(x(u, t), y(u,t)\bigr)$, where the parameter $u$ is taken from a fixed interval $[a,b]$, then we have
\begin{align*}
  \frac{d} {dt} \int_{\cA_i} ydx
  &= \int_{\cA_i} y_t dx + y dx_t\\
  &= \bigl[yx_t\bigr]_{\pd\cA_i} + \int_{\cA_i} y_t dx - x_t dy\\
  &= \bigl[yx_t\bigr]_{\pd\cA_i} + \int_{\cA_i} \bigl(y_t x_s - x_t y_s\bigr) ds\\
\end{align*}
Since the arcs $\cA_i$ evolve by Curve Shortening, we have
\begin{align*}
  x_t &= x_{ss} +\lambda(u, t) x_s \\
  y_t &= y_{ss} +\lambda(u, t) y_s 
\end{align*}
where $\lambda$ is the tangential velocity of the parametrization.  Thus we get
\begin{align*}
  \frac{d} {dt} \int_{\cA_i} ydx
  &= \bigl[yx_t\bigr]_{\pd\cA_i} + \int_{\cA_i} \bigl(y_{ss} x_s - x_{ss} y_s\bigr) ds.
\end{align*}
Finally, if $\theta$ is the tangent angle to the arc, then
\[
  x_s = \cos\theta, \quad y_s = \sin\theta, \quad x_{ss} = -(\sin\theta) \, \theta_s,\quad y_{ss} = (\cos\theta)\, \theta_s.
\]
Hence
\begin{align*}
  \frac{d} {dt} \int_{\cA_i} ydx
  &= \bigl[yx_t\bigr]_{\pd\cA_i} +
    \int_{\cA_i} \bigl( y_{ss}x_s - x_{ss}y_s \bigr) ds \\
  &= \bigl[yx_t\bigr]_{\pd\cA_i} +
    \int_{\cA_i} \theta_s\, ds \\
  &= \bigl[yx_t + \theta\bigr]_{\pd\cA_i}.
\end{align*}
If we now add this over all smooth arcs $\cA_i(t)$, then the $yx_t$ terms cancel because each corner point $A_k(t)$ contributes the same term twice, once from each arc ending at the corner point.  At the end of the asymptotes $y$ is bounded (converges to either $a_0$ or $a_n$), while $x_t$ vanishes.

The $[\theta]_{\pd\cA_i(t)}$ terms together add up to the difference between the changes of $\theta$ along $\cC_j(t)$ and $\Cb(t)$, so that we recover the result in \eqref{eq-area-rate-of-change}.

\bibliography{bibfile.bib}{} \bibliographystyle{plain}

\end{document}